\documentclass[11pt]{amsart}
\usepackage{amsmath,amsthm,amsfonts,amssymb,amscd}
\usepackage{hyperref}

\theoremstyle{plain}
\newtheorem{Main}{Theorem}
\newtheorem{Cor}[Main]{Corollary}

\newtheorem{Thm}{Theorem}[section]
\newtheorem{Lem}[Thm]{Lemma}
\newtheorem{Prop}[Thm]{Proposition}

\newtheorem{Add}[Thm]{Addendum}

\newtheorem{Rk}[Thm]{Remark}

\theoremstyle{definition}
\newtheorem{Defn}[Thm]{Definition}

\newenvironment{pf}{ \begin{proof} }{ \end{proof} }

\DeclareMathAlphabet\EuScript{U}{eus}{m}{n}
\SetMathAlphabet\EuScript{bold}{U}{eus}{b}{n}

\DeclareFontFamily{U}{eus}{\skewchar\font'60}%
\DeclareFontShape{U}{eus}{m}{n}{<-6>eusm5<6-8>eusm7<8->eusm10}{}%
\DeclareFontShape{U}{eus}{b}{n}{<-6>eusb5<6-8>eusb7<8->eusb10}{}%

\newcommand{\Z}{\mathbb{Z}}

\newcommand{\R}{\mathbb{R}}
\newcommand{\C}{\mathbb{C}}

\DeclareMathOperator*{\U}{\mathrm{U}}

\DeclareMathOperator{\Hom}{\mathrm{Hom}}
\DeclareMathOperator{\sym}{Sym}
\DeclareMathOperator{\Pic}{Pic}
\DeclareMathOperator{\vor}{\EuScript{V}}
\DeclareMathOperator{\dbar}{\bar{\partial}}

\DeclareMathOperator{\aut}{Aut}
\DeclareMathOperator{\diff}{Diff}

\DeclareMathOperator{\supp}{Supp}
\DeclareMathOperator{\spinc}{\mathrm{Spin}^c}

\newcommand{\ii}{\mathrm{i}}
\newcommand{\id}{\mathrm{id}}
\newcommand{\Tv}{T^{\mathrm{v}}}
\newcommand{\Th}{T^{\mathrm{h}}}
\newcommand{\univ}{\mathrm{univ}}

\begin{document}

\title{Symplectic fibrations and the abelian vortex equations} 
\author{T. Perutz}
\email{T.Perutz@dpmms.cam.ac.uk}
\date{submitted 8 June 2006; accepted 7 August 2007}
\address{DPMMS, Centre for Mathematical Sciences, University of Cambridge, Wilberforce Road, Cambridge CB3 0WB, United Kingdom.}

\begin{abstract}
The $n$th symmetric product of a Riemann surface carries a natural family of K\"ahler forms, arising from its interpretation as a moduli space of abelian vortices. We give a new proof of a formula of Manton--Nasir \cite{MN} for the cohomology classes of these forms. Further, we show how these ideas generalise to families of Riemann surfaces.

These results help to clarify a conjecture of D. Salamon \cite{Sal} on the relationship between  Seiberg--Witten theory on 3--manifolds fibred over the circle and symplectic Floer homology.
\end{abstract}
\maketitle

\section{Introduction}
\subsection{Relative symmetric products}
Consider a pair of smooth, oriented manifolds $X$ and $S$ with $\dim(X)-\dim(S)=2$, and a proper submersion $\pi\colon X\to S$. Thus $\pi$ is a smooth fibre bundle, and its typical fibre is a compact orientable surface $\Sigma$. 
\begin{Defn}
The $r$th symmetric product bundle, or {\bf relative symmetric product}, $\pi_r\colon \sym^r_S(X)\to S $, is defined to be the quotient by the symmetric group $S_r$ of the fibre product 
\[X^{\times r}_S = \{(x_1,\dots,x_r)\in X^{\times r}: \pi(x_1)=\dots = \pi(x_r)\}\]
with its natural projection to $S$. 
\end{Defn}
\subsubsection{Smooth structures} $\sym^r_S(X)$ is a topological manifold, but it does not inherit a smooth structure from $X$. To make  $\sym^r_S(X)$ a smooth manifold one should choose a complex structure $j$ on the vertical tangent bundle $\Tv X=\ker(D\pi)\subset TX$, compatibly with the orientations. Then the fibres become \emph{complex} manifolds. The smooth atlas on the relative symmetric product is generated by charts which are obtained by fibrewise application of the elementary symmetric functions to `restricted charts' 
\[\Psi\colon D^2\times U \to X.\] 
This means that there is a chart $\psi\colon U\to S$ such that (i) $\pi\circ \Psi =\psi\circ \mathrm{pr}_2$, and (ii) $\Psi\colon D^2\times \{s\}\to X_{\psi(s)}$ is a holomorphic embedding, for each $s\in U$. As observed by Donaldson and Smith \cite{DS}, the existence of such charts is a consequence of the parametrised Riemann mapping theorem. We will write $\sym^r_S(X;j)$ when we want to emphasise that this is the smooth structure being considered.

Different choices, say $j_0$, $j_1$, give distinct smooth structures. However, $\sym^r_S(X;j_0)$
is diffeomorphic to $\sym^r_S(X;j_1)$, as one can see by considering the relative symmetric product of $X\times [0,1]\to S\times [0,1]$, equipped with an interpolating family $j_t$.

\subsubsection{K\"ahler forms} 
The symmetric product $\sym^r(\Sigma)$ of a Riemann surface equipped with a K\"ahler form $\omega$ is itself a K\"ahler manifold. To be precise, a K\"ahler form is determined by a
hermitian line bundle $(L,|\cdot|)$ of degree $r$ over $\Sigma$, together with a real parameter $\tau > 2\pi r/ \int_\Sigma{\omega}$. The reason is that the symmetric product can be identified canonically with a moduli space of abelian vortices, and this has a natural quotient symplectic structure. 

There is a generalisation of this to the case of relative symmetric products. We first fix our conventions concerning families of symplectic manifolds:

\begin{Defn}
(a) A {\bf symplectic fibration} with typical fibre $(M,\omega)$ is a smooth fibre bundle $p\colon X\to S$ together with a vertical two-form $\widetilde{\omega}$, i.e. a section of $\Lambda^2 (\Tv X )^*$, such that each fibre $(X_s,\widetilde{\omega}|X_s)$ is a symplectic manifold isomorphic to $(M,\omega)$.

(b) A {\bf locally Hamiltonian fibration} (LHF) is a triple $(X,p,\Omega)$, where $p \colon X \to S$ is a smooth fibre bundle and $\Omega$ a closed two-form on $X$ such that $(X_s,\Omega| X_s)$ is a symplectic manifold for each $s\in S$.\footnote{The term `locally Hamiltonian fibration' is used in \cite{MS2} in a slightly more restrictive way than here; there it assumed that the base is 2-dimensional and that the form satisfies the normalisation condition introduced by Guillemin and Sternberg.}

\end{Defn}

Relative symmetric products of symplectic surface-fibrations are again symplectic fibrations: if $(p\colon X \to S, \widetilde{\omega})$ is a symplectic fibration with typical fibre $(\Sigma,\omega)$, and one specifies a hermitian line bundle over $X$ of fibrewise degree $r$ and a real parameter, then $\sym^r_S(X)\to S$ becomes a symplectic fibration. \emph{In this paper we show how to promote this functor to locally Hamiltonian fibrations, using the abelian vortex equations.} In doing so we extend Salamon's work \cite{Sal} which applies to bundles over $S^1$. Our method enables one to determine the cohomology classes of the closed forms which arise, in terms of natural operations relating the cohomologies of $X$ and $\sym^r_S(X)$. 

\subsection{Statement of results}
There is a sequence of natural operations sending cohomology classes on $X$ to cohomology classes on the relative symmetric product $\sym^r_S(X)$ of $X\to S$ to classes on $X$. These come about via the universal (or tautological) divisor
\[  \Delta^{\univ} =  {\sym}^r_S(X)\times_S X, \] 
i.e., the locus of pairs $(D,x)$ where $x\in \supp(D)$. This carries a codimension-two homology class relative to boundary, and dually, a cohomology class
\[\delta\in H^2(\sym^r_S(X) \times_S X;\Z).\] 
For example, when $X\to S$ is a holomorphic fibration, $\delta=c_1(\EuScript{O}(\Delta^{\univ}))$. Using the projection maps 
\[\begin{CD}  
 \sym^r_S(X) @<{p_1}<< \sym^r_S(X)\times_S X @>{p_2}>>  X\\ 
 \end{CD}\] 
and cup products in cohomology, define, for each $k\geq 0$, the map
\begin{equation} \label{symm coh ops}
H^*(X;\Z)\to H^{*+2k-2}({\sym}^r_S(X);\Z);\quad c\mapsto c^{[k]} := p_{1!}\big ( (p_2^*c)\smile  \delta^k \big). 
\end{equation}
These operations evidently behave in a natural way under base-change (i.e. pulling back by $S'\to S$).  It is known \cite[Lemma 2.1.1]{Pe1} that
\[ c_1(\Tv {\sym}^r_S(X)) = \frac{1}{2}\left(c_1(\Tv X)^{[1]}+ 1^{[2]}\right).  \]

\begin{Main}\label{Kahler summary}
Let $(X,\pi,\Omega)$ be a proper, locally Hamiltonian surface-fibration over a manifold $S$, and $r$ a positive integer. Choose 
\begin{itemize}
\item
an $\Omega$-positive complex structure $j$ on $\Tv X$;
\item
a hermitian line bundle $(L,|\cdot|)$ over $X$ such that $L|X_s$ has degree $r$ for each $s\in S$, and a unitary connection $A_{\mathrm{ref}}$ on $L$;
\item
a real parameter $\tau$. We require $\tau >2\pi r  a^{-1}$, where $a$ is the symplectic area of a fibre.
\end{itemize}
There is a procedure which associates to these data a closed two-form $v(\Omega,\tau, L)$ on the relative symmetric product $\sym^r_S(X;j)$ which makes it a locally Hamiltonian fibration. This procedure is compatible with restriction of the base $S$. The form
$v(\Omega,\tau,L)$ restricts on each fibre $\sym^r(X_s)$ to the canonical  
K\"ahler form arising from the abelian vortex equations with parameter $\tau$.
Its cohomology class is
\[  [v(\Omega,\tau,L) ] = 2\pi \left( \tau [\Omega]^{[1]} - \pi \,  1^{[2]} \right) \in H^2({\sym}^r_S(X);\R). \]
In particular, the class $[v (\Omega,\tau,L)] $ does not depend on the line bundle $L$. 
\end{Main}

By applying the theorem to fibrations $X\times U \to S\times U$, one sees that there is smooth dependence on parameters. It can be verified without difficulty that, when the base $S$ is the circle, the form $v(\Omega,\tau,L)$ coincides with the one found by Salamon in \cite{Sal}.

\begin{Rk}
In the case where the base $S$ is a point, the result specialises to a formula for the cohomology class of the canonical K\"ahler form on the vortex moduli space
(Theorem \ref{Kahler class}). This formula is due to Manton and Nasir \cite{MN}).  Note, though, that their work relies on a local expansion of the K\"ahler form \cite{Sam} whose derivation has not received the thoroughgoing analytic treatment a pure mathematician would ask for.
\end{Rk}

Some further remarks on the nature of the theorem are in order. The interesting thing is not the existence of closed, fibrewise-K\"ahler two-forms in the specified cohomology class. Indeed, a patching procedure due to Thurston, standard in symplectic geometry, gives an easy construction of such forms. The point is rather that, among such forms, there are some which have a definite geometric (specifically, gauge-theoretic) origin.  This geometric construction is closely related to the Seiberg--Witten equations on fibred 3-- and 4--manifolds: see \cite{Sal} and our discussion of Floer homology below.

Let us say that locally Hamiltonian structures $\Omega_0$, $\Omega_1$ on the same fibre bundle $\pi\colon X\to S$ are {\bf isotopic} if there exists a locally Hamiltonian structure $\Omega\in \Omega^2_{[0,1]\times X}$ on $\pi\times \id \colon [0,1]\times  X\to [0,1]\times S$ with $\Omega|{\{i\}\times X} = \Omega_i$ for $i=0$, $1$. We call LHFs {\bf equivalent} if they are related under the equivalence relation generated by isotopy and two-form-preserving bundle isomorphism.

\begin{Cor} \label{equiv fibrations}
Fix a proper surface-bundle $\pi\colon X\to S$. Choose two sets of data
\[ (\Omega_0, j_0,L_0,|\cdot|_0,A_{\mathrm{ref},0},\tau_0),\quad (\Omega_1, j_1, L_1,|\cdot|_1,A_{\mathrm{ref},1},\tau_1)  \]
as above, and suppose that
$ [\tau _0 \Omega_0]=[\tau_1 \Omega_1]\in H^2(X;\R)$. Then the LHFs
\[(\sym^r_S(X; j_0), \pi_r, v(\Omega,L_0,\tau_0)),\quad (\sym^r_S(X; j_1),\pi_r, v(\Omega,L_1,\tau)) \] are equivalent.   
\end{Cor}

\begin{pf}
Because $\tau_0\Omega_0$ and $\tau_1 \Omega_1$ represent the same cohomology class, the  locally Hamiltonian fibrations $(X,\pi,\tau_0 \Omega_0)$ and $(X,\pi, \tau_1 \Omega_1)$ are isotopic: an isotopy is given by the form $\tau_0 \Omega_0 + d(t\beta )\in Z^2(X\times [0,1])$, where
$\tau_1 \Omega_1 - \tau_0 \Omega_0 = d\beta$. This restricts to the slice $X\times \{t\}$ as $(1-t)\tau_0\Omega_0 + t\tau_1\Omega_1$, and hence is positive on the fibres $X_{s,t}$ of $X\times [0,1]\to S\times [0,1]$. We can give $X\times[0,1]$ a vertical complex structure $J$ by choosing a path $ j_t$ between the given ones.

In the case that $L_1= L_0$, there is a hermitian line bundle $(L,|\cdot|)$ with connection over $X\times [0,1]$ which restricts to $(L_i,|\cdot|_i)$ on the ends. The form 
\[v(\tau_0\Omega_0 + d(t\beta) ,L, 1 ) \] 
on $\sym^r_{S\times[0,1]}(X\times [0,1] ; J)$ restricts on the ends to $v(\tau_i \Omega_i,L_i,1)$. Hence $(\sym^r_S(X; j_0), \pi_r, v(\Omega, L_0,\tau_0))$ is equivalent to $(\sym^r_S(X; j_1), \pi_r, v(\Omega, L_1,\tau_1)$. 

It remains to show that changing the line bundle does not affect things, and for this we may assume that $\Omega_0=\Omega_1$ (write $\Omega$ for this single form) and $j_0 = j_1$. By the theorem, we can write $v(\Omega, L_0,\tau)-v(\Omega,L_1,\tau)=d \gamma$. Then, since $v(\Omega, L_0,\tau)$ and $v(\Omega, L_1,\tau)$ are both K\"ahler, the form $v(\Omega, L_0,\tau) + d(t\gamma)$ on $\sym^r_S(X)\times [0,1]$ gives an isotopy between them.
\end{pf}

\subsection{Floer homology for fibred three-manifolds}
Floer homology for symplectic automorphisms works as follows. Let $\Lambda$ be the universal $\Z/2$--Novikov ring: the ring of formal `series' $\sum_{\lambda\in \R}{ a(\lambda) t^\lambda }$, where $a\colon \R\to \Z/2$ is a function such that $ (-\infty, c] \cap \supp (a)$ is finite for any $c\in \R$. Let $(X,p,\Omega)$ be a LHF over a compact one-manifold $S$, and suppose that its fibres are compact, `weakly monotone'  symplectic manifolds (i.e. $c_1(X_s)$ is positively proportional to the symplectic class, or else $c_1(X_s)$ vanishes on $\pi_2(X_s)$, or else every $S \in \pi_2(X_s)$ has absolute Chern number $| \langle c_1(X_s),[S]\rangle |\geq \dim(X_s)/2-2$). One can then associate with $(X,p,\Omega)$ a $\Lambda$-module 
\[HF_*(X,p; \Omega).\] 
The underlying chain group is freely generated by the set of sections of $X$ which are horizontal for the natural connection determined by $\Omega$ (more precisely, by some  generic perturbation $\Omega'$ of $\Omega$). The differential involves moduli spaces of pseudo-holomorphic sections of $X\times \R \to S\times \R$.

Isotopic LHFs have isomorphic Floer homologies. Two-form-preserving bundle isomorphisms also give isomorphisms in Floer homology. 

This theory is in a sense too rich: different local Hamiltonian structures may give different modules. In  the case where the fibres $X_s$ are complex manifolds, one way to make the theory more manageable is to consider closed two-forms $\Omega$ which are not just fibrewise-symplectic, but actually fibrewise-K\"ahler.  If one also fixes the cohomology class of these forms
then the set of possible choices is a convex set, and $HF_*(X,p;\Omega)$ is independent of the specific choice of $\Omega$. 

An example of this method of making Floer homology manageable occurs in work of  Seidel \cite{Sei}, who applies it to mapping tori of automorphisms of a surface $\Sigma$ of genus $g \geq 2$.  He thereby constructs invariants 
of mapping classes,
\[ \pi_0 \diff^+(\Sigma) \owns [\phi] \mapsto  HF_*( [\phi] ).  \]
Though not well understood, these invariants are far from trivial: Seidel shows that  the identity mapping class $[\id]$ is characterised by the property that, under a natural action by the homology of the surface on Floer homology, $H^2(\Sigma;\Z/2)$ does not annihilate the whole module. 

One way to generalise Seidel's set-up is as follows. Let $\pi\colon Y\to S^1$ be a three-manifold fibred over $S^1$, and consider its relative symmetric product $\sym^r_{S^1}(Y;j)$. We make it an LHF using a closed, fibrewise-K\"ahler two-form drawn from a particular cohomology class. The output will then depend only on the cohomology class chosen: for each class $W$ which restricts to a K\"ahler class on the fibre $\sym^r(\Sigma)$, we get a module $HF_*(\sym^r_{S^1}(Y), \pi_r; W)$. The requirement that the fibres should be weakly monotone forces us to \emph{exclude} the range $g/2\leq r < g-1$, where $g$ is the genus of the fibre.

Let us take $W$ to be one of the classes occurring in our theorem: $W=W(w) = 2\pi \left( w^{[1]} - \pi \,  1^{[2]} \right)$ where $w=[\tau\Omega]$. In this way we obtain a Floer homology module
\[  HF_*(Y,\pi,r; w):= HF_*(\sym^r_{S^1}(Y),\pi_r; W)  \]
by giving only $(Y,\pi,r)$ together with a class $w\in H^2(Y;\R)$ which integrates positively over the fibres of $\pi$. Corollary \ref{equiv fibrations} implies that these modules are well-defined, up to canonical isomorphism.  For discussion of the dependence on $w$ we refer to \cite{Ush}.

Now, we can of course represent $W$ by one of the forms $v(\Omega,\tau,L)$  supplied by the theorem. Doing so is not of any great help in computing Floer homology, but it is highly relevant when we try to understand the relation between the symplectic Floer theory just discussed and the \emph{monopole Floer homology} of the three-manifold $Y$. 

Monopole Floer homology is the Floer theory arising from the Chern--Simon--Dirac functional over a 3--manifold with $\spinc$--structure, a functional whose critical points are precisely the Seiberg--Witten monopoles. Specifically, the name refers to the theory constructed by Kronheimer--Mrowka in their authoritative forthcoming book \cite{KM}. It is a `perturbed' version of monopole Floer homology which is of interest and for this we can again use $\Lambda$ as coefficient ring. 

Salamon's proposal from \cite{Sal}, based on an adiabatic limit computation, 
is that there should be an isomorphism between symplectic and monopole Floer homologies.  Expressing the conjecture in terms of Kronheimer--Mrowka's conventions (and in terms of the notions of this paper) requires a little care because Salamon's conventions differ in various (inessential) ways. If I have accounted correctly for these discrepancies, the statement is that there is an isomorphism between symplectic Floer homology for $\sym^r_{S^1}(Y)$ with the form $v(\Omega,\tau,L)$ (i.e., $HF_*(Y,\pi,r; w)$, where $w=[\tau\Omega]$) and a certain summand in the $\Lambda$-module
\[ HM_*(Y;  - 4 \pi w - 32 \pi^2 c_1(\Tv Y)),\] 
the monopole Floer homology with perturbation class $-4  \pi w - 32 \pi^2 c_1(\Tv Y)$. (The perturbation class is non-zero, providing we assume $g\geq 0$ or $\tau \gg 0$, so that there are no reducible monopoles and only one version of monopole Floer theory.) The summand in question is the direct sum of submodules $HM_*(Y,\mathfrak{t};w)$ where $\mathfrak{t}$ ranges over those $\spinc$-structures on $Y$ for which $\langle c_1(\mathfrak{t}), [\Sigma]\rangle = \chi(\Sigma) + 2r$. 

Let us tie up this discussion. On one hand, we can use pure symplectic geometry to build a group $HF_*(Y,\pi, r; w)$. Specifying the cohomology class on a relative symmetric product---namely, $w^{[1]} - \pi 1^{[2]}$ or a multiple of it---is an essential part of the construction. On the other hand, the existence of the special forms $v(\Omega,\tau,L)$, and Salamon's adiabatic limit, suggest that these modules should have a gauge-theoretic interpretation.

We note finally that the modules $ HF_*(Y,\pi,r; w)$  fit into a field theory for Lefschetz fibrations over surfaces with boundary, which has been studied by M. Usher \cite{Ush} and the author \cite{Pe2} (the latter extends the framework to a larger class of singular fibrations). This too is thought to be intimately related to Seiberg--Witten theory.

\subsection{Acknowledgements}
The work presented here formed a part of my doctoral thesis. I am grateful to my Ph.D. supervisor, Simon Donaldson, for his ideas and advice. Thanks also to Michael Thaddeus for pointing out the Duistermaat--Heckman method, and to Michael Usher for telling me about his related work \cite{Ush}. I acknowledge support from EPSRC Research Grant EP/C535995/1.

\newpage

\section{The vortex equations}
\subsection{Review of moduli spaces of vortices}
Fix a closed Riemann surface $(\Sigma,j)$, a K\"ahler form $\omega\in \Omega^{1,1}_\Sigma$, and a hermitian line bundle $(L,| \cdot |)$ over $\Sigma$, of degree $r>0$. Let $\EuScript{A}(L,|\cdot |$), or $\EuScript{A}(L)$, denote the space of $\U(1)$-connections (an affine space modelled on the imaginary one-forms $\ii \Omega^1_\Sigma$). The gauge group, of smooth maps from $\Sigma$ to $\U(1)$, is denoted by $\EuScript{G}$. Its Lie algebra is $\ii \Omega^0_\Sigma$. The pairing
\[ \Omega^0_\Sigma \otimes \ii \Omega^0_\Sigma \to \R,\quad f\otimes\ii g\mapsto \int_\Sigma{fg \,\omega}\]
embeds $\Omega^0_\Sigma$ into the dual of the $\ii \Omega^0_\Sigma$.
We consider moment maps for Hamiltonian $\EuScript{G}$-actions as maps into  $\Omega^0_\Sigma$.

Connections, sections and gauge transformations are by default $C^\infty$, and the spaces are given their $C^\infty$ topologies. We also need $\EuScript{A}^2_1$, the space of $\U(1)$-connections of Sobolev class  $L^2_1$ (i.e. differing from a smooth one by an $L^2_1$ form); the space of sections $L^2_1(L)$; and the Sobolev gauge group $\EuScript{G}^2_2=L^2_2(\Sigma, \U(1))$. Note that a  map $\Sigma\to \C$ of class $L^2_2$ is continuous, by the Sobolev embedding theorem, and hence has a pointwise norm.

\subsubsection{Action of the gauge group}

The conformal structure $j$ induces a K\"ahler structure on the space of connections $\EuScript{A}(L)$. Its two-form is
\begin{equation}\label{conn 2form}      
(a_1,a_2)\mapsto  \int_\Sigma{\ii a_1\wedge \ii a_2}, \quad a_1,a_2\in \ii\Omega^1_\Sigma.
\end{equation}
The complex structure is the Hodge star $a\mapsto *_j a$. The action of the gauge group $\EuScript{G}$ on $\EuScript{A}(L)$ is Hamiltonian, with (equivariant) moment map
\begin{equation}        
\EuScript{A}(L)\to \ii\Omega^0_\Sigma; \quad A\mapsto * \ii F_A. 
\end{equation}
The symplectic form $\omega$ induces a K\"ahler structure on $\Omega^0_\Sigma(L)$, with two-form
\begin{equation}        
(\phi_1,\phi_2)\mapsto \int_\Sigma{\mathrm{Im}\langle \phi_1,\phi_2 \rangle\, \omega},\quad \phi_1,\phi_2\in \Omega^0_\Sigma(L)       
\end{equation}
and complex structure $\phi\mapsto \ii\phi$. The gauge-action on $\Omega^0_\Sigma(L)$ is Hamiltonian with moment map
\begin{equation}        
\Omega^0_\Sigma(L)\to  \Omega_\Sigma^0; \quad \psi\mapsto \frac{1}{2} |\psi|^2.
\end{equation}
The manifold  $ \EuScript{C}(L):=\EuScript{A}(L)\times\Omega^0_\Sigma(L)$ carries the product K\"ahler structure $\sigma$, which depends on both $j$ and $\omega$. The moment map $m$ for the diagonal $\EuScript{G}$-action is the sum of the moment maps of the factors,
\begin{equation} 
m\colon \EuScript{C}(L)\to \Omega^0_\Sigma,\quad m(A,\psi)=  * i F_A+\frac{1}{2}|\psi|^2 .
\end{equation}
The Chern--Weil formula gives some basic information about this moment map:
\[ \frac{1}{2\pi }\int_\Sigma\tau\omega   
\begin{cases}
< r: & m^{-1}(\tau)=\emptyset; \\
=  r: & m(A,\psi) =\tau \Rightarrow \psi \equiv 0;\\
>  r: & m(A,\psi) =\tau \Rightarrow \psi \not \equiv 0.
\end{cases}\]
In fact, $m$ is submersive at $(A,\psi)$ precisely when $\psi \not\equiv 0$, which is also the locus on which the gauge-action is free. When $\int_\Sigma\tau\omega >2\pi r$, the free gauge-action on $\mu^{-1}(\tau)$ admits local slices (see below), so the K\"ahler quotient $m^{-1}(\tau)/ \EuScript{G}$ is a K\"ahler manifold.

\subsubsection{The vortex equations}
The {\bf vortex equations} with parameter $\tau$ are the following coupled equations for a pair $(A,\psi)\in\EuScript{C}(L)$:
\begin{align}
\label{v1} \dbar_A\psi&=0       && \text{in  } \Omega^{0,1}_\Sigma(L), \\ 
\label{v2} m(A,\psi) &=\tau   && \text{in  } \Omega^0_\Sigma.
\end{align}
Individually, we will refer to them as the Cauchy--Riemann equation and the moment map equation. The space of solutions $\widetilde{\vor}(L,\tau)$ is invariant under $\EuScript{G}$, and the quotient space \[\vor(L,\tau) := \widetilde{\vor}(L,\tau)/\EuScript{G}\] is called the vortex moduli space.
The fundamental results about $\vor(L,\tau) $ are as follows. 
\begin{Prop}
Assume that $ \int{\tau\omega}>2\pi r $.
\begin{enumerate}
\item[(a)]
The space $\vor(L,\tau)$ is a finite-dimensional, complex---therefore smooth and K\"ahler---submanifold of $m^{-1}(\tau)/\EuScript{G}$.
\item[(b)]
The map
\[      Z\colon \vor(L,\tau) \to \sym^r(\Sigma), \quad  [A,\psi]\mapsto \psi^{-1}(0)    \] 
is an isomorphism of complex manifolds.
\end{enumerate}
\end{Prop}
The unitary connection $A$ induces a holomorphic structure on $L$: a local section is
holomorphic if and only if it lies in $\ker \dbar_A$. By means of the holomorphic structure, one attaches  multiplicities to points of $\psi^{-1}(0)$, so that $\psi$ has $r$ zeros in all. This makes sense of $Z$. We write $L_A$ for $L$ with this holomorphic structure.

Item (a) is proved by an elliptic regularity argument, and we shall say a little more about it. As for (b), the statement that $Z$ is bijective is an existence and uniqueness theorem for solutions to the vortex equations. This is the heart of the theorem, and various proofs are known, see e.g. Jaffe and Taubes \cite{JT}, Garc\'ia-Prada \cite{GP}.

The `degenerate' case, where $\int{\tau \omega}=2\pi r$, is also interesting:
\begin{Add}
When $ \int{\tau\omega}=2\pi r $, the moduli space
\[\vor(L,\tau)=\{(A,0):\ii F_A =  \tau\omega \}/\EuScript{G}\]
is a finite-dimensional, complex---therefore smooth and K\"ahler---submanifold of $m^{-1}(\tau)/\EuScript{G}$. The map
\[ \vor(L,\tau) \to \Pic^L(\Sigma);\quad [A,0]\mapsto L_A \]
is an isomorphism of complex manifolds.
\end{Add}
Here $\Pic^L(\Sigma)$ is the Picard torus of holomorphic structures on $L$.
\subsubsection{Smoothness of the moduli space}
This is a standard application of elliptic theory. We run through it briefly in preparation for the family version considered later; see \cite{Sal} for some more details. 

The tangent space to the affine space $\EuScript{C}^2_1$ is the space of pairs $(a,\phi)$, where $a$ is an imaginary one-form, $\phi$ a section, both of class $L^2_1$. One obtains local slices for the action of $\EuScript{G}$ by imposing the Coulomb gauge condition
\begin{equation}\label{gauge} 
d^*(\ii a)+ \mathrm{Im}\langle \psi,\phi \rangle =0, 
\end{equation}
which says that $(a,\phi)$ is orthogonal to the gauge-orbit of $(A,\psi)$. Note that the left-hand side is gauge-equivariant. The linearisations of the two vortex equations at the solution $(A,\psi)$ are 
\begin{equation}\label{linearised}
\dbar_A\phi + a^{0,1}\psi =0, \quad *\ii da +\mathrm{Re}\langle \psi,\phi\rangle=0.
\end{equation}
The second of these and (\ref{gauge}) are real and imaginary parts of the single equation
\begin{equation}\label{gauss}
 \dbar^* (a^{0,1}) - \frac{1}{2}\langle\psi,\phi \rangle=0. 
 \end{equation}
Hence the space of solutions to equations (\ref{linearised}, \ref{gauge}) is the kernel of the $\C$-linear differential operator
\begin{equation} 
D_{(A,\psi)}: \; (a,\phi)\mapsto (\dbar_A\phi + a^{0,1}\psi,  \dbar^* (a^{0,1}) - \frac{1}{2}\langle\psi,\phi \rangle). 
\end{equation}
Now, $D_{(A,\psi)}$ is a compact perturbation of the Fredholm operator \
\[(a,\phi)\mapsto (\dbar_A\phi ,\dbar^*(a^{0,1}) ),\] 
which has index $(r+1-g)- (1-g)=r$ (over $\C$). Hence  $D_{(A,\psi)}$ is also Fredholm of index $r$. It is surjective (this can be seen by computing $D_{(A,\psi)}^* D_{(A,\psi)}$, see \cite{Sal}), so its kernel has constant rank $r$. 

From this point it is straightforward to check, using the implicit function theorem, that $\vor(L,\tau)$ is a differentiable submanifold of $m^{-1}(\tau)/\EuScript{G}$. Since its tangent spaces
$\ker(D_{A,\psi})$ are complex linear, it is a complex submanifold.

\subsection{The K\"ahler class on the vortex moduli space}
As we have seen, the moduli space $\vor(L,\tau)$ is a complex manifold equipped with a canonical K\"ahler form $\sigma_\tau$. We write $\sigma_\tau$ also for its pullback by $Z^{-1}$, a K\"ahler form on $\sym^r(\Sigma)$. The target of this section is to determine its cohomology class. 

A $(2-p)$-cycle $\zeta$ in $\Sigma$ gives rise to a closed subset $\delta_\zeta\subset \sym^r(\Sigma)$ representing a $(2r-p)$-cycle: $\delta_\zeta$ consists of divisors $D\in \sym^r(\Sigma)$ such that $\mathrm{mult}_{x}(D) = \mathrm{mult}_x(\zeta)$ for all $x\in \Sigma$. Using this map followed by Poincar\'e duality on $\sym^r(\Sigma)$, we obtain a map 
\[\nu_p \colon H_p(\Sigma;\Z) \to H^{2-p}(\sym^r(\Sigma);\Z) .\] 
It is well-known that $\nu_1$ is an isomorphism. When $p=2$, an isomorphism
 \[   H_0(\Sigma;\Z) \oplus \Lambda^2 H_1(\Sigma;\Z) \stackrel{\cong}{\to} H^2(\sym^r(\Sigma);\Z)  \] 
is given by
\[ (a , b\wedge c )\mapsto \nu_0(a) +   \nu_1(b) \cup \nu_1(c). \]
We define
\begin{itemize}
\item
$\eta \in H^2(\sym^r(\Sigma);\Z) $ to be the class corresponding to the point class in $H_0(\Sigma;\Z)$;
\item
$\theta\in H^2(\sym^r(\Sigma);\Z)$ to be the class corresponding to the cup-product form on $H^1(\Sigma;\Z)$ (here we think of the cup-product form as an element of $\Hom(\Lambda^2 H^1(\Sigma;\Z),\Z)=\Lambda^2 H_1(\Sigma;\Z)$). 
\end{itemize}
Often we conflate these integral classes with their images in real cohomology.

\begin{Main}\label{Kahler class}
The equation
\[\frac{1}{2\pi} [\sigma_\tau] = \left( \int_\Sigma{\tau \omega}\right) \eta + 2\pi (\theta-r\eta)\]
holds in $H^2(\sym^r(\Sigma);\R)$.
\end{Main}
As already mentioned, this formula was found by Manton--Nasir \cite{MN}.
Our (quite different) method of proof is to exhibit connections on two line bundles over the orbit space of irreducible pairs, $\EuScript{C}^*/\EuScript{G}$. The Chern classes of these line bundles restrict to $\eta$ and $\theta-r\eta$ on $\vor(L,\tau)$, while the appropriate linear combination of their curvature forms restricts exactly to the form $\sigma_\tau$. 

\subsubsection{Cohomology of the orbit space}
We write $\EuScript{C}^*=\EuScript{C}^*(L)$ for the space of pairs $(A,\psi)\in\EuScript{C}(L)$ with $\psi$ not identically zero, $\EuScript{B}^*$ for the orbit space $\EuScript{C}^*/\EuScript{G}$, and $i\colon \vor(L,\tau) \to \EuScript{B}^*$ for the inclusion.
 \begin{Lem}
 $i$ induces a surjection on cohomology, and an isomorphism on $H^{\leq 2}$.
 \end{Lem}
\begin{pf}
Using the cohomology slant product operation, define
\begin{align*}
\mu_{\EuScript{B}}\colon &H_*(\Sigma;\Z)\to H^{2-*}(\EuScript{B}^*;\Z), && h\mapsto c_1(\EuScript{L}_{\EuScript{B}})/h,\\
\mu_{\sym}\colon &H_*(\Sigma;\Z)\to H^{2-*}(\sym^r(\Sigma);\Z), && h\mapsto c_1(\EuScript{L}_{\sym})/h.
\end{align*}
Here the line bundle $\EuScript{L}_{\EuScript{B}}\to \EuScript{B}^*\times \Sigma$ is  $\EuScript{L}_{\EuScript{B}}=\widetilde{\EuScript{L}_{\EuScript{B}}}/\EuScript{G}$, where the equivariant line bundle $\widetilde{\EuScript{L}_{\EuScript{B}}}\to \EuScript{C}^*\times \Sigma$ is the pullback of $L\to \Sigma$; and  $\EuScript{L}_{\sym}\to \sym^r(\Sigma)\times \Sigma$ is the topological line bundle corresponding to the universal divisor $\Delta^{\univ}\subset \sym^r(\Sigma)\times \Sigma$. 

These maps extend uniquely to ring homomorphisms 
\begin{align*}  
&  \Lambda^*H_1(\Sigma;\Z)\otimes_{\Z}\Z[H_0(\Sigma)] \to H^*(\EuScript{B}^*;\Z),\\
&  \Lambda^*H_1(\Sigma;\Z)\otimes_{\Z}\Z[H_0(\Sigma)] \to H^*(\sym^r(\Sigma);\Z),  \end{align*}
since the ring on the left is freely generated by $H_0(\Sigma;\Z)\oplus H_1(\Sigma;\Z)$. These are homomorphisms of \emph{graded} rings where the grading on the left is characterised by the property that $ H_i(\Sigma;\Z)$ has degree $2-i$. The first of these two maps is an isomorphism \cite[pp. 539--545]{AB}. The second is surjective, since the image of $\mu_{\sym}$ contains $H^1(\sym^r(\Sigma);\Z)$ and the class $\eta$, and these generate the cohomology ring. 

To prove the lemma it suffices to show that $i^* \circ \mu_1=Z^*\circ \mu_{\sym}$. This follows from the fact that $(i\times 1)^*\EuScript{L}_{\EuScript{B}}$ is isomorphic to $(Z\times 1)^*\EuScript{L}_{\sym}$. To see that these bundles are isomorphic, observe that the former has a tautological section which vanishes precisely along $\Delta^{\univ}$.
\end{pf}

It is convenient to have some notation to hand for integral (co)homology classes on $\Sigma$. 
Let $e_0 \in  H_0(\Sigma)$ be the class of a point,
$e_2 \in H_2(\Sigma)$ the orientation class. Let $e^0 \in H^0(\Sigma)$, $e^2\in
H^2(\Sigma)$ be their duals. Let 
$\{\alpha_i,\beta_j\}_{1\leq i,j \leq g} $
be a symplectic basis for $H_1(\Sigma)$, and $\{ \alpha^i,\beta^j\}_{1\leq i,j \leq g}$ the dual basis for $H^1(\Sigma)$.

Now put
\begin{equation} 
\widetilde{\eta}=\mu_1(e_0),\quad \widetilde{\theta} = \sum_{i=1}^g{\mu_1(\alpha_i)\cup \mu_1(\beta_i)}. 
\end{equation}
 \begin{Lem}
 $c_1(\EuScript{L}_{\EuScript{B}})^2/e_2=2 r\widetilde{\eta}-2 \widetilde{\theta}$ in $H^2(\EuScript{B}^*;\Z)$.
 \end{Lem}
\begin{pf}
The group $H^2(\EuScript{B}^*\times \Sigma;\Z)$ is the direct sum of its K\"unneth components $H^0(\EuScript{B}^*;\Z)\otimes H^2(\Sigma;\Z)$, $H^1(\EuScript{B}^*;\Z)\otimes H^1(\Sigma;\Z)$ and $H^2(\EuScript{B}^*;\Z)\otimes H^0(\Sigma;\Z)$. The Chern class $c_1(\EuScript{L}_{\EuScript{B}})$ is \emph{tautologically} the sum of

 \begin{align*}
 \mu_1(e_2) \otimes e^2  &\in  H^0(\EuScript{B}^*;\Z)\otimes H^2(\Sigma;\Z),\\ 
  \sum_{i=1}^g{(\mu_1(\alpha_i)\otimes \alpha^i + \mu_1(\beta_i) \otimes \beta^i)}  
                &\in  H^1(\EuScript{B}^*;\Z)\otimes H^1(\Sigma;\Z),\\
 \mu_1(e_0) \otimes e^0 &\in  H^2(\EuScript{B}^*;\Z)\otimes H^0(\Sigma;\Z).
 \end{align*}
Let us call these terms $A$, $B$ and $C$ respectively. Note that $A = r.1 \otimes
 e^2$ (by definition of $\EuScript{L}_{\EuScript{B}}$) and $C= \widetilde{\eta}\otimes
 e^0$. 
The K\"unneth isomorphism is compatible with cup products, providing that
one uses the \emph{graded} tensor product of graded rings. Thus
$A \cup C  = r \widetilde{\eta} \otimes e^2 = C\cup A $, and 
\[B^2=  \left(\sum_{i=1}^g {\mu_1(\alpha_i)\cup \mu_1(\beta_i)-\mu_1(\beta_i) \cup \mu_1(\alpha_i)}
\right)\otimes e^2 = -2 \widetilde{\theta} \otimes e^2.\]
Hence $c_1(\EuScript{L}_{\EuScript{B}})^2/e_2=2(r\widetilde{\eta}-\widetilde{\theta})$.
\end{pf}
\subsubsection{A connection on \texorpdfstring{$\EuScript{L}_{\EuScript{B}}$}{the universal line bundle}}

We now write down a canonical connection $\nabla$ on $\EuScript{L}_{\EuScript{B}}$, and compute its curvature. This calculation is modelled on that of Donaldson and Kronheimer \cite[p. 195]{DK}. We will use the curvature form, together with its wedge-square, to construct a closed two-form on $\EuScript{B}^*$, representing a known cohomology class, whose restriction to $\vor(L,\tau)$ is $\sigma_\tau$. 

The connection $\nabla$ is concocted from two ingredients:
\begin{itemize}
\item a certain unitary, $\EuScript{G}$-invariant connection $\widehat{\nabla}$ on the line bundle $\mathrm{pr}_2^*L \to \EuScript{A}^*\times \Sigma$;
\item
a certain connection $\Gamma$ on the principal $\EuScript{G}$-bundle $\EuScript{C}^*\to \EuScript{B}^*$, pulled back to $\EuScript{B}^*\times\Sigma$. 
\end{itemize}
As explained in \cite{DK}, such data determine a connection $\nabla$ on the quotient line bundle $\EuScript{L}_{\EuScript{B}}\to \EuScript{B}^*\times\Sigma$, characterised by the condition
\[ (\nabla _v s) \, \hat{}  = \widehat{\nabla}_{\hat{v}} (\hat{s})  \]
for local sections $s$ and vector fields $v$, where $\hat{\cdot}$ denotes $\Gamma$-horizontal lifting.

The connection $\widehat{\nabla}$ is trivial in the $\EuScript{C}^*$-directions and tautological in the $\Sigma$-directions. To amplify: a section of $\mathrm{pr}_2^*L$ is a map $s\colon \EuScript{C}^*\times \Sigma\to L$ with $s(A,\psi,x)\in L_x$, and at the point $(A,\psi,x)$,
 \begin{equation} \label{taut conn}
 \widehat{\nabla}_{(a,\phi,v)}s = d_{A,v}(s|{\{A, \psi\}\times \Sigma})(x)+ \left(\left.\frac{d}{dt}\right|_{t=0}s(A+ta,\psi+t\phi,x)\right)(x). 
 \end{equation}
The curvature of $\widehat{\nabla}$ is given by
\begin{align} 
&F_{\widehat{\nabla}}((0,0,u),(0,0,v))   =F_A(u,v),  \notag \\
&F_{\widehat{\nabla}}((a,\phi,0),(0,0,v))  = \langle a ,v\rangle, \\
&F_{\widehat{\nabla}}((a,\phi,0),(a',\phi',0)) = 0.\notag
\end{align}
We can obtain a connection on $\EuScript{C}^*\to\EuScript{B}^*$ from our gauge-fixing condition: the horizontal space over $[A,\psi]$ is the kernel of the linear operator 
\[(a,\phi)\mapsto d^*(\ii a) -\mathrm{Im}\langle\psi,\phi\rangle.\] 
To write down the connection one-form $\Gamma$, we need the Green's operator $G_\psi$ associated to the Laplacian
\[ \Delta_\psi = d^* d + |\psi|^2 \colon \quad \Omega^0_\Sigma \to \Omega^0_\Sigma.\] 
$\Delta_\psi$ is surjective (since $d^*d$ maps onto the functions of mean-value zero), inducing an isomorphism of $\ker(\Delta_\psi)^\perp$ with $\Omega^0_\Sigma$; its inverse is $G_\psi$.
\begin{Lem}
The connection one-form $\Gamma$ is given by
\[ \Gamma_{(A,\psi)}(a,\phi) = \ii  G_\psi (d^*\ii a - \mathrm{Im}\langle\psi,\phi\rangle)  
\in \ii \Omega^0_\Sigma. \]
\end{Lem}
\begin{pf}
This form has the correct kernel, so to justify the assertion one simply observes that it is invariant under $\EuScript{G}$:
\[ \Gamma_{(A,\psi,x)}( - d f , f\psi,0 )= f , \quad f\in \ii \Omega^0_\Sigma. \]
\end{pf}
In accordance with the general pattern explained in \cite{DK}, the curvature of the quotient connection $\nabla$ on $\EuScript{L}_{\EuScript{B}}\to \EuScript{B}^*\times \Sigma$ is given by
\begin{align} 
&F_{\nabla}((0,0,u),(0,0,v))   =F_A(u,v), \notag  \\
&F_{\nabla}((a,\phi,0),(0,0,v))  = \langle a ,v\rangle, \label{quot curv} \\
&F_{\nabla}((a_1,\phi_1,0),(a_2,\phi_2,0)) = 2 \ii G_\psi(d^*\ii b-\mathrm{Im}\langle \psi ,\chi\rangle).\notag
\end{align}
Here $(a_1,\phi_1)$ and $(a_2,\phi_2)$ are vector fields on $\EuScript{B}^*$ which are \emph{horizontal} with respect to $\Gamma$; their Lie bracket is $(b,\chi)$. 

\begin{Lem}
Suppose $(a_1,\phi_1)$ and $(a_2,\phi_2)$ are horizontal. Then \[d^*(\ii b)-\mathrm{Im}\langle\psi,\chi \rangle =- \mathrm{Im}\langle\phi_1,\phi_2\rangle.\]
\end{Lem}
\begin{pf}
Denote the pair $(A+t a_1,\psi+t\phi_1)$ by $c_t$. Then, at $(A,\psi)$,
\begin{align*} 
& b = \frac{1}{t}\big( a_2(c_0)-a_2(c_t)\big)+o(t) , \\ 
& \chi = \frac{1}{t}\big( \phi_2(c_0)-\phi_2(c_t) \big ) +o(t) 
\end{align*}
as $t\to 0$. But at $c_t$, $d^*\,\ii a_2=\mathrm{Im}\langle \psi+t\phi_1, \phi_2(c_t) \rangle$, and from this one obtains
\[d^*( \ii b )- \mathrm{Im}\langle\psi,\chi\rangle = - \lim_{t\to 0}\, \mathrm{Im} \langle \phi_1(c_0),\phi_2(c_t) \rangle = - \mathrm{Im} \langle \phi_1(c_0),\phi_2(c_0)\rangle.  \]
\end{pf} 
\subsubsection{Two-forms as curvature integrals}

We are now in a position to write down closed two-forms representing $c_1(\EuScript{L}_{\EuScript{B}})/e_0$ and $c_1(\EuScript{L}_{\EuScript{B}})^2/e_2$ in de Rham cohomology.

Note: \emph{in this paragraph we insist that the tangent vectors $(a_j,\phi_j)$ are horizontal.}
 
The class $c_1(\EuScript{L}_{\EuScript{B}})$ has the Chern--Weil representative $\ii F_\nabla/2\pi$, so  
\begin{equation} 
c_1(\EuScript{L}_{\EuScript{B}})/e_0  = \left [  \frac{1}{2\pi} \int_\Sigma \ii F_\nabla \wedge \omega_0\right], \quad \text{where }\int_\Sigma{\omega_0}=1.
\end{equation} 
Explicitly, this representative for $c_1(\EuScript{L}_{\EuScript{B}})/e_0$ is the two-form
\begin{equation} 
((a_1,\phi_1),(a_2,\phi_2))\mapsto  \frac{1}{\pi [\omega]}\int_\Sigma {G_\psi(\mathrm{Im}\langle \phi_1,\phi_2 \rangle)\,\omega}. 
\end{equation}  
Similarly,
\begin{equation} 
 c_1(\EuScript{L}_{\EuScript{B}})^2/e_2 = \left[\frac{1}{4\pi^2}\int_\Sigma{\ii F_\nabla\wedge \ii F_\nabla}\right]. \end{equation} 
This integral involves the product of the first and third curvature terms, and the square of the second term. So $c_1(\EuScript{L}_{\EuScript{B}})^2/e_2$ has the representative 
\begin{equation}  
((a_1,\phi_1),(a_2,\phi_2))\mapsto  \frac{1}{\pi^2} \int_\Sigma{G_\psi(\mathrm{Im}\langle \phi_1,\phi_2 \rangle) \ii F_A } -\frac{1}{2\pi^2}\int_\Sigma{\ii a_1\wedge \ii a_2}.
\end{equation}
Notice the appearance of an expression familiar from (\ref{conn 2form}) as the second term. 

At this point we impose the moment map equation, restricting these forms and classes to the locus where $m(A,\psi)=\tau$. On that locus, the class
\begin{equation} 4 \pi^2 (\widetilde{\theta} - r\widetilde{\eta}) + 2\pi [\tau \omega]\widetilde{\eta}  =  2 \pi   \big(-\pi c_1(\EuScript{L}_{\EuScript{B}})^2/e_2 + [ \tau \omega] c_1(\EuScript{L}_{\EuScript{B}})/e_0\big) 
\end{equation}  
is represented by the form
\begin{align}
& \int_\Sigma{\ii a_1\wedge \ii a_2} + 2\int_\Sigma{G_\psi(\mathrm{Im}\langle \phi_1,\phi_2 \rangle)\big( \tau \omega - \ii F_A \big) } \notag \\
 = &\int_\Sigma{\ii a_1\wedge \ii a_2} + \int_\Sigma{G_\psi(\mathrm{Im}\langle \phi_1,\phi_2 \rangle)|\psi|^2 \omega }\notag \\
=& \int_\Sigma{\ii a_1\wedge \ii a_2} +  \int_\Sigma{\mathrm{Im}\langle \phi_1,\phi_2 \rangle \omega}\notag \\
=& \sigma( (a_1,\phi_1), (a_2,\phi_2)) \label{manipulation}. 
\end{align}
(Recall that $\sigma$ is our standard K\"ahler form on $\EuScript{C}^*$). The penultimate equality uses the observation that, because the Laplacian
of a function $f$ has mean value zero,
\[ \int_\Sigma{f \omega}= \int_\Sigma{\Delta_\psi G_\psi (f)\,\omega}= \int_\Sigma {|\psi|^2G_\psi (f)}\,\omega. \]

\begin{pf}[Proof of Theorem \ref{Kahler class}]
What we have just found is that the class $ 2\pi([\tau \omega]\widetilde{\eta} + 2\pi (\widetilde{\theta}-r \widetilde{\eta} ))$ on $\EuScript{B}^*$, restricted to $m^{-1}(\tau)/\EuScript{G}$, is equal to $[\sigma_\tau]$. Restricting further to the vortex moduli space, we find that the class of our preferred K\"ahler form is
\[ 2\pi ([ \tau\omega]\widetilde{\eta} + 2\pi (\widetilde{\theta}-r \widetilde{\eta} )) |{\vor(L,\tau)}\in H^2(\vor(L,\tau);\R). \]
Hence, pulling back by $Z$, we find that the class of our K\"ahler form on $\sym^r(\Sigma)$ is $2\pi([ \tau \omega]\eta + 2\pi(\theta-r \eta))$, which is the formula we have been working towards.
\end{pf}

\subsection{The Duistermaat--Heckman formula}
The Duistermaat--Heckman formula \cite{DH} for the variation of cohomology of symplectic quotients gives another proof that the cohomology class $[\sigma_\tau]$ varies linearly with $\tau$---provided that $\tau$ is a constant function---and computes the slope.

Suppose that one has a Hamiltonian $S^1$-action on $(M,\omega)$, with moment map $\mu\colon M\to \mathfrak{t}^*$. Here $\mathfrak{t}=\mathrm{Lie}(S^1)$. Identify $\mathfrak{t}^*$ with $\R$ so that the lattice dual to $\exp^{-1}(1)\subset\mathfrak{t}$ corresponds to $\Z\subset \R$.
Suppose that $\mu$ is proper, and that its restriction to $\mu^{-1}(\ell)$ is submersive, for some open interval $\ell\subset \R$. The family of symplectic quotients $(M_t,\omega_t)_{t\in\ell}$ is then a trivial fibre bundle, and a trivialisation gives an identification of the cohomology of $M_t$ with that of a fixed fibre $M_s$. The identification is canonical, hence $\{[\omega_t]\}_{t\in\ell}$ can be considered as a family of classes on $M_s$. Suppose that $S^1$ acts freely on $\mu^{-1}(s)$, so that $\mu^{-1}(s)\to M_s$ is a principal circle-bundle, with Chern class $c\in H^2(M_s;\R)$. The Duistermaat--Heckman formula says that 
\begin{equation}\label{DH} 
\frac{d}{dt}[\omega_t] =  2\pi c. 
\end{equation}
We apply this with $M=\bigcup_{\tau\in\R}{\widetilde{\vor}(L,\tau)/\EuScript{G}_0}$, where $\EuScript{G}_0$ is the based gauge group $\{ u\colon \Sigma\to \U(1): u(x)=1 \}$, $x\in \Sigma$ an arbitrary basepoint, and $\tau\in \R$ stands for a constant function on $\Sigma$. The circle acts by constant gauge transformations. We take $\ell=(2\pi r/\int{\omega},\infty)$; the Chern class $c$ of $\widetilde{ \vor} (L,\tau)/\EuScript{G}_0\to \widetilde{ \vor} (L,\tau)/\EuScript{G}$ is $\eta$. Formula \ref{DH} gives
\[ \frac{d}{d\tau}[\sigma_\tau] = 2\pi \left(\int_\Sigma{\omega}\right) \eta, \]
which is consistent with our result. One can formally recover the constant term $4\pi^2 (\theta-r \eta )$ by specialising to the degenerate parameter $\tau = 2\pi r/\int{\omega}$ (for which the formula $[\sigma_\tau]=4\pi^2\theta$ is easily verified); however, justifying this formal manipulation would need further thought.

Since Duistermaat and Heckman's proof identifies the variation in the symplectic forms with the curvature of a connection on $\mu^{-1}(s)\to M_s$, the two methods are perhaps not so different as they appear.

\section{Families of vortex moduli spaces}
\subsection{Construction of the vortex fibration}\label{constr}
(a) Suppose that $X\to S$ is a smooth fibre bundle, where $X$ and $S$ are connected and oriented, and that the typical fibre is a compact surface $\Sigma$. Let $L\to X$ be a principal $\U(1)$-bundle, and assume that $L|{X_s}\to X_s$ has degree $r>0$.

Consider $X\to S$ as a fibration with structure group $\diff^+(\Sigma)$. Putting $P=L|{X_s}$, we can consider the composite map $L\to X\to S$ as a fibration with typical fibre $P$ and structure group $\diff^+_P(\Sigma)$. The latter is the group of pairs $(\tilde{g},g)$, where $\tilde{g}\in \aut(P)$ is an automorphism covering $g\in \diff^+(\Sigma)$, so it is an extension of $\diff^+(\Sigma)$ by the gauge group.

There are natural left actions of $\diff^+_P(\Sigma)$ on the space of connections $\EuScript{A}(P)$ and on the space of sections $\Omega^0_\Sigma(P)$. These arise through the covariance of connections and of sections; representing a connection by its one-form $A\in\Omega^1_P$, we have 
\[ \tilde{g}.A =\tilde{g}^{-1*}A; \quad  \tilde{g}.\psi=\tilde{g}\circ\psi\circ g^{-1}. \]
One can then form the associated fibrations  
\[ L\times_{\diff^+_P(\Sigma)} \EuScript{A}(P)\to S,\quad L\times_{\diff^+_P(\Sigma)} \Omega^0_\Sigma(P)\to S,  \]
with structure group $\diff^+_P(\Sigma)$. These may be thought of as the bundles of connections (resp. sections) along the fibres of $X\to S$:
\begin{align*} 
L\times_{\diff^+_P(\Sigma)} \EuScript{A}(P) & \cong\{(s,A): s\in S,\,A\in\EuScript{A}_{X_s}(L_s)\},\\
L\times_{\diff^+_P(\Sigma)} \Omega^0_\Sigma(P) & \cong 
\{ (s,\psi): s\in S,\, \psi \in\Omega^0_{X_s}(L_s) \}.
\end{align*}
The first of these has the special property that it is a \emph{symplectic} fibration: its structure group is reduced to the symplectic automorphism group of $\EuScript{A}(P)$.

Other fibrations can be derived from these basic ones. The space 
\[ \EuScript{C}(P) = \EuScript{A}(P)\times \Omega^0_\Sigma(P\times_{U(1)}\C) , \]
comprising pairs $(A,\psi)$ where $\psi$ is a section of the line bundle associated with $P$, is also a $\diff^+_P(\Sigma)$-space (the action is the diagonal one), and so is
\[ \EuScript{B}(P) = \EuScript{C}(P)/\EuScript{G}, \]
because $\EuScript{G}$ acts on $\EuScript{C}(P)$ as a subgroup of $\diff^+_P(\Sigma)$.
The associated fibrations are
\begin{align*}  \EuScript{C}_{X/S}(L) & :=  L\times_{\diff^+_P(\Sigma)}  \EuScript{C}(P),\\
\EuScript{B}_{X/S}(L) & :=  L\times_{\diff^+_P(\Sigma)}  \EuScript{B}(P).
\end{align*}

(b) Suppose now that $X\to S$ is itself a symplectic fibration, i.e. that its structure group is reduced to  $\aut(\Sigma,\omega)$ for some area form $\omega$. Then the structure group $L\to S$ is reduced to $\aut_P(\Sigma,\omega)$, the group of pairs $(\tilde{f},f)$ with $f^*\omega=\omega$, and $\EuScript{C}_{X/S}(L)\to S$ is again a symplectic fibration. Note that $P\times_{\U(1)}\C$ is a \emph{hermitian} line bundle, so our formula for the symplectic form on $\EuScript{C}(P)$ makes sense.

Let $\{j_s\in\EuScript{J}(X_s,\omega_s)\}_{s\in S}$ be a smooth family of complex structures, compatible with the symplectic forms. The moment map $m \colon \EuScript{C}(P)\to \ii\Omega^2_\Sigma$, $(A,\psi)\mapsto * \ii F_A + |\psi|^2 /2$, generalises to a bundle map over $S$,
\[  m_{X/S} \colon  \EuScript{C}_{X/S}(L)\to L\times_{\diff^+_P(\Sigma)}\Omega^0_\Sigma. \]
We now take $\tau$ to be a \emph{constant}. Then we have a sub-bundle
\[  m_{X/S}^{-1}(\tau) \subset \EuScript{C}_{X/S}(L), \]
projecting to a sub-bundle  $\pi (m_{X/S}^{-1}(\tau)) \subset \EuScript{B}_{X/S}(L)$ under the quotient map $\pi\colon \EuScript{C}_{X/S}(L)\to\EuScript{B}_{X/S}(L)$, and $\pi (m_{X/S}^{-1}(\tau\omega))\to S$ has structure group $\aut_P(\Sigma,\omega)$.

(c) We now impose a fibred version of the Cauchy--Riemann equation. This differs from what we have done so far in that it cannot be expressed in terms of associated bundles.

The total space of the {\bf vortex fibration} $\vor_{X/S}(L,\tau) \to S$ is the space of triples $[s,A,\psi] \in \pi (m_{X/S}^{-1}(\tau)) $ satisfying $\dbar_{j_s,A}\psi=0$. It maps to $S$ in the obvious way. The fibre over $s$ can be identified with the vortex moduli space $\vor_{X_s}(L|{X_s},\tau)$, and so with $\sym^r(X_s)$.
\begin{Lem}
The space $\vor_{X/S}(L,\tau)$ has a structure of smooth manifold which makes the projection $p\colon \vor_{X/S}(L,\tau)\to S$ a smooth submersion, hence a fibre bundle. 
\end{Lem} 
\begin{pf}
The linearisation of the defining equations for $\vor_{X/S}(L,\tau)$, and the fibrewise gauge-fixing condition, define an $\R$-linear operator $ D_{(s,A,\psi)}$:
\begin{equation}
D_{(s,A,\psi)}(v,a,\phi) = D_{A,\psi}(a,\phi)+ P(v),\quad v\in T_s S. 
\end{equation}
Here $P$ is the $0$th-order operator $P(v)=\frac{1}{2}\ii(d_A\psi)\circ \frac{\partial j}{\partial v}$. The operator $D_{(s,A,\psi)}$ is thus Fredholm, of real index $2r+\dim(S)$, and surjective (since $D_{(A,\psi)}$ is). The kernel of $D_{(s,A,\psi)}$ is the putative tangent space to $\vor_{X/S}(L,\tau)$ at $(s,A,\psi)$, and the projection $\pi \colon \ker D_{(s,A,\psi)}\to T_s S $ is putatively the derivative of $p$. Note that $\pi$ is surjective, because its kernel is exactly $\ker D_{A,\psi}$, which we know has dimension $2r$.  Now the standard elliptic theory which we sketched above here gives smoothness of the vortex fibration and of the map $p$.
\end{pf}

\subsection{Line bundles and cohomology operations.}

Let $ \widetilde{\EuScript{L}}_{\EuScript{B}}\to \EuScript{C}^*_{X/S}\times_S X $ be the pullback of the line bundle $L\to X$. It is an equivariant line bundle under the fibrewise gauge-action, and so descends to a line bundle 
\[ \EuScript{L}_{\EuScript{B}}\to \EuScript{B}^*_{X/S} \times_S X. \] 
The universal divisor $ \Delta^{\univ} \subset {\sym}^r_S(X)\times_S X$
corresponds to a unique line bundle 
\[\EuScript{L}_{\sym}\to \sym^r_S(X)\times_S X.\] 

\begin{Lem}
There is a natural isomorphism
\[ (i\times 1)^*\EuScript{L}_{\EuScript{B}} \to Z^* \EuScript{L}_{\sym}, \] 
where $i$ is the inclusion of $\vor_{X/S}(L,\tau)$ in $\EuScript{B}^*_{X/S}$,
and $Z$ the natural isomorphism of $\vor_{X/S}(L,\tau)$ with $\sym^r_S{X}$.
\end{Lem}

\begin{pf}
The section 
\[ ( [A,\psi],x) \mapsto [\psi(x)] \]
of $(i\times 1)^*\EuScript{L}_{\EuScript{B}} $ vanishes precisely along $Z^{-1}(\Delta)$.
\end{pf}
Using these two line bundles one can construct operations
\begin{align*} 
& H^*(X)\to H^{*+2k-2} ( \EuScript{B}^*_{X/S}(L)) , &&  c\mapsto \tilde{c}^{[k]}, \\
& H^*(X)\to H^{*+2k-2} ({\sym}^r_{S}(X)), && c\mapsto c^{[k]}.
\end{align*}
defined for arbitrary coefficient rings. The second of these earlier was discussed earlier (Equation \ref{symm coh ops}). Introduce the projections 
\[\begin{CD} 
 \EuScript{B}^*_{X/S}  @<{p_1}<< \EuScript{B}^*_{X/S} \times_S X @>{p_2}>> X,\\
  \sym^r_S(X)  @<{p_1}<< \sym^r_S(X)\times_S X @>{p_2}>> X,
\end{CD} \]
and set
\begin{align} 
\tilde{c}^{[k]} = p_{1!}(c_1(\EuScript{L}_{\EuScript{B}})^{k}\cup p_2^* c), \\
c^{[k]} = p_{1!}(c_1(\EuScript{L}_{\sym})^{k}\cup p_2^* c). 
\end{align}
Because of the relation between $\EuScript{L}_{\EuScript{B}}$ and $\EuScript{L}_{\sym}$, we have $i^*\tilde{c}^{[k]} = Z^* c^{[k]}.  $

\subsection{Associated fibrations as locally Hamiltonian fibrations}

In Section \ref{constr}, we constructed various associated fibrations within the category of \emph{symplectic} fibrations---fibre bundles with symplectic forms on the fibres. Our next task is to refine these constructions to the category of \emph{locally Hamiltonian} fibrations. The vortex fibration will then become a LHF by restricting a closed two-form defined on a larger space. The cleanest way that I have found to do this is to `reverse-engineer' our cohomology calculation for the vortex moduli space. This goes as follows.

We need a fibrewise-equivariant connection $\widehat{\nabla}$ on the bundle 
\[\widetilde{\EuScript{L}}_{\EuScript{B}}\to\EuScript{C}^*_{X/S}\times_S X. \]
To obtain one, choose a connection $A_{\mathrm{ref}}$ on $L\to X$. We define $\widehat{\nabla}$ to be the unique connection which restricts to the natural one (\ref{taut conn}) on each fibre over $S$, and which is given by $A_{\mathrm{ref}}$ on $\Th X$.

In conjunction with the fibrewise gauge-fixing condition, $\widehat{\nabla}$ defines a quotient connection $\nabla$ on $\EuScript{L}_{\EuScript{B}}\to \EuScript{B}^*_{X/S}\times_S X$. 
\begin{Defn}
We define the closed two-form $\tilde{v}(\tau \Omega,L)$ on $\EuScript{B}^*_{X/S}$ by
\begin{equation}  
\tilde{v}(\tau \Omega,L) = 2\pi \int_{X/S}{\mathrm{i}F_\nabla \wedge ( \tau \Omega- \frac{\pi}{2} \mathrm{i}F_\nabla)}
\end{equation}
We define $v(\Omega,\tau,L)$ to be the restriction of $i^*\tilde{v}(\tau\Omega,L)$ to the vortex fibration $\vor_{X/S}(L,\tau)$. 
\end{Defn}
Let us clarify the integration symbol here. Projection on the first factor makes $\EuScript{B}^*_{X/S}\times_S X$ a fibre bundle over $\EuScript{B}^*_{X/S}$. The fibre over a point of $\EuScript{B}^*_{X/S}\times_S X$ which lies over $s\in S$ is $X_s$. It therefore makes sense to integrate down the fibres of $\EuScript{B}^*_{X/S}\times_S X\to\EuScript{B}^*_{X/S}$. In particular, a closed four-form $\alpha$ on $\EuScript{B}^*_{X/S}\times_S X$ gives rise to a closed two-form 
\[\int_{({B}^*_{X/S}\times_S X) /(\EuScript{B}^*_{X/S})}{\alpha}.\]  
We write this more compactly as $\int_{X/S}{\alpha}$.

Bearing in mind that integration along the fibre corresponds to the cohomology push-forward, we can read off the cohomology class of $\tilde{v}(\tau \Omega,L)$:
\[ [\tilde{v}(\tau \Omega,L)] =  2\pi ([\tau\tilde{\Omega}]^{[1]} - \pi \tilde{1}^{[2]}). \]
Forming $\tilde{v}(\tau\Omega,L)$ is obviously compatible with restricting the base $S$. By our earlier calculations, the form $i^*\tilde{v}(\tau\Omega,L)$ on the vortex bundle restricts to the preferred K\"ahler form on each fibre. Thus
\[  [v (\tau \Omega,L)]  =  2\pi  ([\tau\Omega]^{[1]} - \pi 1^{[2]} ). \]
Theorem \ref{Kahler summary} is now an immediate consequence of what we have done.

\end{document}